\documentclass[12pt]{article}

\usepackage{amstext    }
\usepackage{amsthm    }
\usepackage{a4}
\usepackage[mathscr]{eucal}

\usepackage{amsmath}
\usepackage{amssymb}
\usepackage{amscd}

\numberwithin{equation}{section}

\newtheorem{theorem}{Theorem}[section]

\newtheorem{proposition}[theorem]{Proposition}
\newtheorem{corollary}[theorem]{Corollary}
\newtheorem{lemma}[theorem]{Lemma}

\newcommand{\cali}[1]{\mathscr{#1}}

\newcommand{\PGL}{{\rm PGL}}

\newcommand{\supp}{{\rm supp}}

\newcommand{\dist}{{\rm \ \! dist}}

\newcommand{\ddc}{dd^c}

\newcommand{\DSH}{{\rm DSH}}

\newcommand{\id}{{\rm id}}

\newcommand{\Ac}{\cali{A}}

\newcommand{\Cc}{\cali{C}}

\newcommand{\Ec}{\cali{E}}
\newcommand{\Fc}{\cali{F}}

\newcommand{\Hc}{\cali{H}}

\newcommand{\FS}{{\rm FS}}

\newcommand{\C}{\mathbb{C}}

\newcommand{\R}{\mathbb{R}}
\renewcommand{\P}{\mathbb{P}}

\title{Equidistribution speed for endomorphisms of projective spaces}
\author{Tien-Cuong Dinh and Nessim Sibony}

\begin{document}
\maketitle

\begin{abstract}
Let $f$ be a non-invertible holomorphic endomorphism of $\P^k$, $f^n$ its
iterate of order $n$ and $\mu$ the equilibrium measure of $f$. We 
estimate the speed of convergence in the following known result. 
If $a$ is a Zariski generic point in $\P^k$, the
probability measures, equidistributed on the preimages of $a$ under $f^n$,
converge to $\mu$ as $n$ goes to infinity.
\end{abstract}




\section{Introduction}

Let $\P^k$ denote the complex projective space of dimension $k$. Let
$f:\P^k\rightarrow\P^k$ be a holomorphic endomorphism of
algebraic degree $d\geq 2$. The map $f$ defines a ramified covering of
degree $d^k$ over $\P^k$.  Let $\omega_\FS$ denote the Fubini-Study
form on $\P^k$ normalized so that $\omega_\FS^k$ is a probability
measure. Let $f^n:=f\circ \cdots\circ f$ ($n$ times) be the iterate of
order $n$ of $f$. It is well-known that $d^{-kn} (f^n)^*(\omega_\FS^k)$
converge to a probability measure $\mu$ which is totally invariant:
$d^{-k}f^*(\mu)=f_*(\mu)=\mu$. We refer to the surveys \cite{DinhSibony4,
  Sibony} for the basic dynamical properties of such maps.

Consider a point  $a$ in $\P^k$. 
We are interested on the asymptotic distribution of the fibers
$f^{-n}(a)$ of $a$ when $n$ goes to infinity. More
precisely, if $\delta_a$ denotes the Dirac mass at $a$,
$d^{-kn}(f^n)^*(\delta_a)$ is the probability measure
which is equidistributed on the fiber $f^{-n}(a)$. The points in
$f^{-n}(a)$ are
counted with multiplicities. Here is our first result.

\begin{theorem} \label{th_equi_1} 
Let $f$ be a holomorphic endomorphism of algebraic degree $d\geq 2$ of
$\P^k$. Let $\mu$ be the equilibrium measure of $f$ and $1<\lambda<d$
a constant. There is an invariant proper analytic subset $E_\lambda$,
possibly empty, of
$\P^k$ such that if $a$ is a point out of $E_\lambda$ and if 
$\varphi$ is a $\Cc^\alpha$ function on $\P^k$
with $0<\alpha \leq 2$, then
$$|\langle d^{-kn} (f^n)^*(\delta_a)-\mu,\varphi\rangle|\leq A 
\Big[1 +\log^+{1\over \dist(a,E_\lambda)}\Big]^{\alpha/2}\|\varphi\|_{\Cc^\alpha}
\lambda^{-\alpha n/2},$$
where $A>0$ is a constant independent of $n$, $a$, $\varphi$ and $\log^+(\cdot):=\max(0,\log(\cdot))$. 
\end{theorem}

Note that the distance $\dist(a,E_\lambda)$ is with respect to the
Fubini-Study metric on $\P^k$. When $E_\lambda$ is empty, by
convention, this distance is the diameter of $\P^k$ which is a finite
number. A priori, the constant $A$ depends on
$\lambda$ and $\alpha$. Note also that we have the estimate
$$A \Big[1+\log^+{1\over \dist(a,E_\lambda)}\Big]^{\alpha/2}
\leq A \Big[1+\log^+{1\over\dist(a,E_\lambda)}\Big]$$
for $0<\alpha\leq 2$. 

Let $\Hc_d(\P^k)$ denote the set of holomorphic endomorphisms of
algebraic degree $d$ of $\P^k$. We can identify it with a Zariski
open set of a projective space.  It is shown in
\cite[Lemma 5.4.5]{DinhSibony3}, that for $f$ in a dense Zariski open set
$\Hc^\lambda_d(\P^k)$ in $\Hc_d(\P^k)$, the multiplicities of points in the
fibers of $f^{n_0}$ are smaller than $d^{n_0}\lambda^{-n_0}$ for some
$n_0\geq 1$. We will see in Section \ref{section_excep} that in this
case, $E_\lambda$ is empty. The following consequence of Theorem
\ref{th_equi_1} gives a precise version of \cite[prop. 5.4.13]{DinhSibony3}.

\begin{corollary}
Let $1<\lambda<d$ be a constant and $f$ an element of
$\Hc_d^\lambda(\P^k)$. Then for every
$a$ in $\P^k$ and for $\varphi$ a $\Cc^\alpha$ function on $\P^k$
with $0<\alpha \leq 2$, we have
$$|\langle d^{-kn} (f^n)^*(\delta_a)-\mu,\varphi\rangle|\leq A 
\|\varphi\|_{\Cc^\alpha}
\lambda^{-\alpha n/2},$$
where $A>0$ is a constant independent of $n$, $a$ and $\varphi$.
\end{corollary}

Recall that for $f$ in $\Hc_d(\P^k)$,  there is a maximal proper analytic subset $\Ec$ of $\P^k$
which is totally invariant under $f$, i.e. $f^{-1}(\Ec)=f(\Ec)=\Ec$,
see Section \ref{section_excep} for details. Here is our second result.

\begin{theorem} \label{th_equi_2} 
Let $f$, $\mu$ and $\Ec$ be as above. 
There is a constant $\lambda>1$ such that if $a$ is a point out of $\Ec$ and if 
$\varphi$ is a $\Cc^\alpha$ function on $\P^k$
with $0<\alpha \leq 2$, then
$$|\langle d^{-kn} (f^n)^*(\delta_a)-\mu,\varphi\rangle|\leq A 
\Big[1+\log^+{1\over \dist(a,\Ec)}\Big]^{\alpha/2}\|\varphi\|_{\Cc^\alpha}
\lambda^{-\alpha n/2},$$
where $A>0$ is a constant independent of $n$, $a$ and $\varphi$. In
particular, $d^{-kn} (f^n)^*(\delta_a)$ converge to $\mu$
locally uniformly for $a\in\P^k\setminus\Ec$. 
\end{theorem}

The following result is a direct consequence of Theorem
\ref{th_equi_2}.

\begin{corollary} \label{cor_equi}
Let $f$, $\mu$ and $\Ec$ be as above. Then $d^{-kn} (f^n)^*(\delta_a)$
converge to $\mu$ if and only if $a\not\in\Ec$. 
\end{corollary}

Note that in dimension 1, this result was proved by Brolin
\cite{Brolin} for
polynomials, by Lyubich \cite{Lyubich} and Freire-Lopes-Ma{\~n}{\'e} 
\cite{FLM} for general maps.
A version of our theorems was also proved by Drasin-Okuyama
in \cite{DrasinOkuyama} using Nevanlinna theory.

In higher dimension, a partial result in this direction
was obtained by Forn\ae ss and the second author in
\cite{FornaessSibony1}. It shows the equidistribution property for $a$
outside a set of zero Lebesgue measure (a pluripolar set). 
Corollary \ref{cor_equi} was announced by Briend-Duval in
\cite{BriendDuval1}. Their proof shows  the equidistribution
property for $a$ outside a countable union of hypersurfaces (the orbit
of the critical values of $f$). 
The first complete proof of this corollary was given
by the authors 
in \cite{DinhSibony1}, see also \cite{DinhSibony2}. 

Our proof there is valid
for a much more general setting (polynomial-like maps and maps on
singular varieties) and is separated into two parts. The first
one improves a geometrical method due to Lyubich in dimension 1 and developed by
Briend-Duval \cite{BriendDuval1} in higher dimension. The second part,
quite different from the dimension 1 case,
shows the existence of an analytic exceptional set $\Ec$ and some extra
properties which are useful in the proof. This exceptional set is
still not well understood. Very recently, Briend-Duval
\cite{BriendDuval2} showed
that using a dynamical argument, as in our work \cite{Dinh,
  DinhSibony1, DinhSibony2}, one can obtain a short proof of Corollary
\ref{cor_equi}.

In this paper, we will use known properties of the exceptional set
but we will replace the geometrical method with a
pluripotential one which is stronger.
The strategy was already introduced by Forn\ae ss and the second author in
\cite{FornaessSibony2} for the equidistribution of hypersurfaces. In
\cite{DinhSibony3}, we showed that this strategy can be extended to
the case of varieties of arbitrary dimension, in particular, for the
equidistribution of points that we consider. The novelty here is
that we obtain precise quantitative results.

\section{P.s.h. functions and regularization}

We refer to \cite{Demailly, Hormander} for the basic properties of
plurisubharmonic (p.s.h. for short) functions and the theory of
positive closed currents. Recall that a function $\varphi$
on $\P^k$ is {\it quasi-p.s.h.} if it is locally equal to the difference of
a p.s.h. function and a smooth function. In particular, such a
function is defined everywhere, bounded from
above and is in all $L^p$ spaces with $1\leq p<\infty$. We also have
$\ddc\varphi\geq -c\omega_\FS$ for some constant $c>0$. A Borel set
$E$ in $\P^k$ is {\it pluripolar} if it is contained in the pole set
$\{\varphi=-\infty\}$ of a quasi-p.s.h. function $\varphi$. 

Recall that a function $\varphi$ on $\P^k$, defined out of a
pluripolar set, is {\it d.s.h.} if it is equal to the difference
of two quasi-p.s.h. functions. We identify two d.s.h. functions if
they are equal out of a pluripolar set. 
We summarize here some properties of these functions, see \cite{DinhSibony6} for details.
If $\varphi$ is d.s.h., 
there are two positive closed $(1,1)$-currents $S^\pm$ such that
$\ddc\varphi= S^+-S^-$. The following expression defines a norm on the
space of d.s.h. functions:
$$\|\varphi\|_\DSH:=|\langle\omega_\FS^k,\varphi\rangle|+\inf\|S^\pm\|,$$
where the infimum is taken over $S^\pm$ as above. 

Recall that the mass
of a positive closed $(p,p)$-current $S$ on $\P^k$ is given by
$\|S\|:=\langle S,\omega_\FS^{k-p}\rangle$. It depends only on the
cohomology class of $S$ in $H^{p,p}(\P^k,\C)$. So, the masses of
$S^+$ and $S^-$ are equal since these currents are in the same class. 
Note that $\Cc^2$ functions are d.s.h. and we have
$\|\cdot\|_\DSH\lesssim \|\cdot\|_{\Cc^2}$. 

It is not difficult to show that we can write 
$$\varphi=\varphi^+-\varphi^-+c$$
out of a pluripolar set, where $c$
is a constant and $\varphi^\pm$ are quasi-p.s.h. such
that $\max\varphi^\pm=0$. Moreover, we have the estimates 
$$\ddc\varphi^\pm\geq 
-M\|\varphi\|_\DSH\omega_\FS \quad \mbox{and}\quad |c|\leq
M\|\varphi\|_\DSH$$ 
with some constant $M>0$ independent of $\varphi$. 

We say that a probability measure $\nu$ on $\P^k$ is {\it PB} if
quasi-p.s.h. functions are integrable with respect to $\nu$. Such a
measure $\nu$ has no mass on pluripolar sets and 
d.s.h. functions are also integrable with respect to $\nu$. When $\nu$ is PB, the following norm is equivalent to
$\|\cdot\|_\DSH$,
$$\|\varphi\|_\nu:=|\langle\nu,\varphi\rangle|+\inf\|S^\pm\|.$$
Recall that the equilibrium measure $\mu$ of $f$ is PB.
In fact, it satisfies a stronger property. On the space of d.s.h. functions, consider the following
{\it weak topology}: a sequence of d.s.h. functions $(\varphi_n)$
converges to $\varphi$ if $\|\varphi_n\|_\DSH$ are uniformly bounded
and $\varphi_n\rightarrow\varphi$ in the sense of currents. In that case,
$\varphi_n\rightarrow\varphi$ in all $L^p$ spaces with $1\leq
p<\infty$. The measure $\mu$ is {\it PC},
i.e. $\varphi\mapsto\langle\mu,\varphi\rangle$ is continuous for the
above topology on $\varphi$.  
We refer to \cite{DinhSibony4} for a recent survey on this theory.

A classical {\it exponential estimate}
implies that if $\varphi$ is quasi-p.s.h. such that $\ddc \varphi\geq
-\omega_\FS$ and $\max\varphi=0$ then 
$$\langle \omega_\FS^k,e^{\alpha|\varphi|}\rangle \leq M$$
for some constants $\alpha>0$ and $M>0$ independent of $\varphi$ \cite{Hormander}. It
is not difficult to see that 
the estimate still holds when $\varphi$ belongs to a
bounded family of d.s.h. functions. The constants $M$ and $\alpha$
depend on the family. The following proposition is crucial in
our estimates and can be extended to H{\"o}lder continuous functions.

\begin{proposition} \label{prop_horm}
Let $\Fc$ be a bounded family of d.s.h. functions. Then, there is a constant
$A_0>0$ such that 
$$\|\varphi\|_\infty\leq A_0(1+\log^+\|\nabla\varphi\|_\infty)$$
for every $\Cc^1$ function $\varphi$ in $\Fc$.
\end{proposition}
\proof
Let $A_0>2$ be a constant large enough. If the above estimate were
false, then
there is a function $\varphi$ in $\Fc$ and a point $a$ in $\P^k$ such
that 
$|\varphi(a)|\geq A_0\log\|\nabla\varphi\|_\infty\geq A_0$.
It follows that
$|\varphi|\geq A_0\log\|\nabla\varphi\|_\infty-1\geq {1\over 2}A_0\log\|\nabla\varphi\|_\infty$ on the disc centered
at $a$ of radius
$\|\nabla\varphi\|_\infty^{-1}$. This contradicts the
exponential estimates for d.s.h. functions.
\endproof

We will use the automorphism group $\PGL(\C,k+1)$ of $\P^k$ in order to regularize
functions. This is a complex Lie group
which acts transitively on $\P^k$. Let $u$ denote a
local holomorphic coordinates system of $\PGL(\C,k+1)$ with $\|u\|<1$
such that $u=0$ at identity $\id\in\PGL(\C,k+1)$. 
The automorphism of coordinates $u$ is denoted by $\tau_u$.
Let $\rho(u)$
be a smooth probability measure with compact support in
$\{\|u\|<1\}$. For any $0<\theta<1$ define
$\rho_\theta(u)$ as the image of this measure under $u\mapsto \theta u$. This is also a smooth probability
measure with support in $\{\|u\|<\theta\}$. If $\phi$ is an $L^1$ 
function on $\P^k$, define {\it the $\theta$-regularized function}
$\phi_\theta$ of $\phi$ by
$$\phi_\theta:=\int_{\PGL(\C,k+1)} (\phi\circ \tau_u) d\rho_\theta(u).$$
It is not difficult to show that $\phi_\theta$ is smooth and converge to
$\phi$ in the $L^1$ topology when $\theta$ tends to 0. The
proof of the following lemma is left to the reader, see also
\cite[Prop. 2.1.6]{DinhSibony3}.

\begin{lemma}  \label{lemma_reg}
There is a constant $A_1\geq 1$ independent of $\phi$ and
  $\theta$ such
  that 
$$\|\phi_\theta\|_{\Cc^1} \leq A_1\|\phi\|_{L^1}\theta^{-7k^2}.$$
Moreover, if $\Ac(\phi,x,r):=\sup_y |\phi(y)-\phi(x)|$ for $y$ such
that $\dist(x,y)\leq r$, then
$$|\phi_\theta(x)-\phi(x)|\leq \Ac(\phi,x,\eta \theta),$$
where $\eta>0$ is a constant independent of $\phi$, $x$ and $\theta$.
\end{lemma} 

We will also need the following lemma.

\begin{lemma} \label{lemma_reg_dsh}
Let $\Fc$ denote the family of d.s.h. functions $\varphi$
  such that $\langle\mu,\varphi\rangle=0$ and $\|\varphi\|_\mu\leq 1$. Then
  $\varphi_\theta-\langle\mu,\varphi_\theta\rangle$ is also a function in $\Fc$.
\end{lemma}
\proof
Define $\varphi_u:=\varphi\circ \tau_u$ and
$\varphi'_u:=\varphi_u-\langle\mu,\varphi_u\rangle$. Since
$\varphi_\theta-\langle\mu,\varphi_\theta\rangle$ is an average on
$\varphi'_u$, it is enough to show that $\varphi'_u$ is in $\Fc$. We
use here that $\mu$ is PC. It is
clear that $\langle\mu,\varphi_u'\rangle=0$. Since $\varphi$ is in $\Fc$,
there are positive closed $(1,1)$-currents $S^\pm$ such that
$\ddc\varphi=S^+-S^-$ and $\|S^\pm\|\leq 1$. We have
$\ddc\varphi_u:=\tau_u^*(S^+)-\tau_u^*(S^-)$. On the other hand, since
the automorphisms of $\P^k$ act trivially on the Hodge cohomology
groups of $\P^k$, we have $\|\tau_u^*(S^\pm)\|=\|S^\pm\|\leq 1$. Indeed,
the mass of a positive closed current depends only on its cohomology
class. Therefore, $\|\varphi'_u\|_\mu\leq 1$. 
\endproof

\section{Exceptional sets} \label{section_excep}

Let $\kappa_n(x)$ denote the local topological degree of $f^n$ at $x$
for $n\geq 0$. More precisely, for $z$ generic near $f^n(x)$,
$f^{-n}(z)$ has $\kappa_n(x)$ points near $x$. 
The functions $\kappa_n$ satisfy $\kappa_{n+m}=(\kappa_m\circ
f^n)\kappa_n$. So, they form a multiplicative cocycle.
Define 
$$\kappa_{-n}(x):=\max_{y\in f^{-n}(x)}\kappa_n(y).$$
These functions are upper semi-continuous with respect to the Zariski
topology on $\P^k$. We summarize here some properties obtained in \cite{Dinh, Favre}. 
The sequence $\kappa_n^{1/n}$ converges to a
function $\kappa_+$. The sequence
$\kappa_{-n}^{1/n}$ converges to a
function $\kappa_-$ which is upper semi-continuous with respect to the
Zariski topology. Moreover, for any $\delta>1$, the level set
$\{\kappa_-\geq \delta\}$ is an invariant proper analytic subset of
$\P^k$. Each point in $\{\kappa_+\geq \delta\}$ is sent to
$\{\kappa_-\geq \delta\}$ by some iterate of $f$. Since $\kappa_+\circ
f=\kappa_+$, we have $f^{-1}\{\kappa_+\geq \delta\}\subset
\{\kappa_+\geq\delta\}$ and in general,
$\{\kappa_+\geq\delta\}$ is not an analytic set.

Define 
$$E_\lambda:=\{\kappa_-\geq d\lambda^{-1}\}\quad\mbox{and}\quad
\Omega:=\P^k\setminus E_\lambda.$$
Then, there is a constant $1<\delta_0<d\lambda^{-1}$ such that
$E_\lambda=\{\kappa_-\geq \delta_0\}$. 
Let $n_0\geq 1$ be an integer such that either $\{\kappa_{-n_0}\geq
\delta_0^{n_0}\}\cap\Omega$ is empty or its dimension $q$ is minimal. 
The definition of $E_\lambda$ and Proposition 2.4 in \cite{Dinh} imply that
this analytic set is in fact empty. So,
$\kappa_{-n_0}<\delta_0^{n_0}$ out of $E_\lambda$. We
will use the following lemma for $Y:=E_\lambda$ or $Y:=\Ec$, they are invariant. 

\begin{lemma} \label{lemma_dist}
There is a constant $A_2\geq 1$ such that for every subsets $X$ and $Y$ of
$\P^k$, we have
$$\dist(f^{-n}(X),f^{-n}(Y))\geq A_2^{-n}\dist(X,Y).$$
In particular, if $f(Y)\subset Y$, we have   
$$\dist(f^{-n}(X),Y)\geq A_2^{-n}\dist(X,Y).$$
\end{lemma} 
\proof
It is enough to take $A_2$ large enough so that $\dist(f(x),f(y))\leq
A_2\dist(x,y)$. If $X,Y$ are two subsets of $\P^k$, we have
$\dist(f^{-1}(X),f^{-1}(Y))\geq A_2^{-1}\dist(X,Y)$. We obtain the
result by induction. 
\endproof

The lemma implies that $\dist(f^{-n_0}(a),E_\lambda)\geq
A_2^{-n_0}\dist(a,E_\lambda)$ and an analogous property for $\Ec$. 
So, in order to prove the main results, we can replace $f$ with $f^{n_0}$
and assume that $n_0=1$. Therefore, the local topological degree of $f$ is
smaller than $\delta_0$ at every point in 
$\P^k\setminus f^{-1}(E_\lambda)$. We can also replace $f$
with an iterate $f^n$ of $f$, $\lambda$, $\delta_0$, $d$ with
$\lambda^n$, $\delta_0^n$ and
$d^n$ in order to assume that $20k^2\delta_0<d\lambda^{-1}$.  

Recall that the measure $\mu$ is a wedge-product of positive closed
$(1,1)$-currents with H{\"o}lder continuous local potentials. We easily
deduce from the estimate in \cite[Prop. 2.3.6]{DinhSibony3} the
following property, see also \cite[Prop. 3.3]{DinhNguyenSibony}.

\begin{lemma} \label{lemma_tub}
Let $V_t$ denote the $t$-neighbourhood of $E_\lambda$,
$0<t<1$, i.e. the set of points $x$ such that
$\dist(x,E_\lambda)\leq t$. Then, 
$$\mu(V_t)\leq A_3 t^\beta,$$ 
where $A_3\geq 1$ and
$\beta>0$ are two constants independent of $t$.
\end{lemma}

It is shown in \cite{Dinh,DinhSibony1, DinhSibony2} that there is a maximal
proper analytic set $\Ec$, possibly empty, which is totally invariant
under $f$, i.e.
$f^{-1}(\Ec)=f(\Ec)=\Ec$. The following property is deduced from the
construction of $\Ec$ as the maximal totally invariant analytic set
\cite[Theorem 6.1]{DinhSibony2}. 

\begin{lemma} \label{lemma_excep}
Let $X$ be a proper analytic set of $\P^k$. If $a$ is a
  point in $\P^k\setminus \Ec$, then there is an integer $n\geq 0$ such that
  $f^{-n}(a)\not\subset X$. 
\end{lemma}

\section{Lojasiewicz inequalities}

Let $B_r$ denote the ball of center 0 and of radius $r$ in $\C^k$. Let
$\pi$ denote the canonical projection from
$\C^k\times\C^k$ onto its first factor.  
Recall the following version of the Lojasiewicz inequality, see
\cite{DinhSibony2, FornaessSibony2}.

\begin{proposition} \label{prop_loj_1}
Let $X$ be an analytic subset of $B_1\times B_1$
   of pure dimension $k$ and $s_0$ a fixed integer. Assume that
  $\pi:X\rightarrow B_1$
  defines a ramified covering of degree $s\leq s_0$ over
  $B_1$. Then there is a constant $c>0$ such that if $x,y$ are two
  points in $B_{3/4}$ we can write
$$\pi^{-1}(x)\cap X=\{x_1,\ldots, x_s\}\quad \mbox{and}\quad 
\pi^{-1}(y)\cap X=\{y_1,\ldots, y_s\}$$
with $\|x_i-y_i\|\leq c\|x-y\|^{1/s}$. Moreover, the constant
$c$ depends on $s_0$ but not on $X$. 
\end{proposition}

Note that the points in the fibers $\pi^{-1}(x)\cap X$ and
$\pi^{-1}(y)\cap X$ are repeated according to their
multiplicities. We deduce from this proposition the following result.

\begin{proposition} \label{prop_loj_2}
Let $X$ be an analytic subset of $B_1\times B_1$
   of pure dimension $k$ and $m$ an integer. Assume that
  $\pi:X\rightarrow B_1$ defines a ramified covering of degree $s$ over
  $B_1$. Let $Z\subset B_1$ be a proper analytic set such that the
  multiplicity of every point in $\pi^{-1}(x)\cap X$ is at most
  equal to $m$ for $x \in B_1\setminus Z$. Then there are
  constants $c>0$, $N\geq 1$ such that for any $0<t<1$ and
  all $x,y\in B_{1/2}$ with $\dist(x,Z)\geq t$, $\dist(y,Z)\geq
  t$,  we can write
$$\pi^{-1}(x)\cap X=\{x_1,\ldots, x_s\}\quad \mbox{and}\quad 
\pi^{-1}(y)\cap X=\{y_1,\ldots, y_s\}$$
with $\|x_i-y_i\|\leq ct^{-N}\|x-y\|^{1/m}$.
\end{proposition}

Fix a constant $A\geq 1$ large
enough. By Proposition \ref{prop_loj_1}, if $\widetilde z$ is a point in
$\pi^{-1}(z)\cap X$ with $z\in B_1$ then
$$\dist(\widetilde z,\pi^{-1}(x)\cap X)\leq A\|x-z\|^{1/s}.$$
For $z\in B_1$, write $\pi^{-1}(z)\cap X=\{z_1,\ldots,z_s\}$. Consider
the multi-indices $I=\{i_1,\ldots,i_{m+1}\}\subset \{1,\ldots,s\}$ of
length $m+1$. Define 
$$h_I(z):=\sum_{1\leq p,q\leq m+1} \|z_{i_p}-z_{i_q}\|.$$ 
We can see $h_I(z)$ as a multi-valued function on $z$.

\begin{lemma} \label{lemma_technique}
There is an integer $M\geq 1$ such that for $z\in B_{1/2}$ we have
$$h_I(z)\geq A^{-1} \dist(z,Z)^M.$$
\end{lemma}
\proof
Let $X_1$ denote the analytic set of points $P=(z,z_I)$ in
$B_1\times (\C^k)^{(m+1)^2}$, where the coordinates of
$z_I$ are defined by $z_{i_p}-z_{i_q}$ with $z_i$ and $1\leq
p,q\leq m+1$ as above.  
It can be seen as a ramified covering over $B_1$. 
Define
$$X_2:=B_1\times \{0\}\quad \mbox{and}\quad X_3:=Z\times
(\C^k)^{(m+1)^2}.$$
We have $X_1\cap X_2\subset X_3$.  A Lojasiewicz inequality \cite[p.14
and p.62]{Malgrange} implies that for $z\in B_{1/2}$ and for $A, M$ large enough
$$\dist(P,X_2)\geq A^{-1} \dist(P,X_1\cap X_2)^M \geq A^{-1}
\dist(P,X_3)^M = A^{-1} \dist(z,Z)^M.$$
We infer 
$$h_I(z)\geq \dist(P,X_2)\geq A^{-1} \dist(z,Z)^M.$$
This completes the proof of the lemma.
\endproof

\noindent
{\bf End of the proof of Proposition \ref{prop_loj_2}.}
Define $N:=Ms$ and fix a constant $\gamma>0$
small enough. We can assume that 
$\|x-y\|\leq {1\over 2}\gamma t^{Nm}$. Otherwise, the lemma is
obvious. Let $B$ denote the ball of center
$x$ and of radius $\gamma t^{Nm}$ in $\C^k$. We will show that each
connected component of $X\cap\pi^{-1}(B)$ defines a ramified covering
of degree $\leq m$ over $B$.
Then it is enough to apply
Proposition \ref{prop_loj_1} to this component. In order to apply
that proposition, we have to use a coordinate dilation on $B$ which leads to
the factor $t^{-N}$ in the estimate $\|x_i-y_i\|\leq
ct^{-N}\|x-y\|^{1/m}$.  

Fix a point $a=(x,x')$ in $\pi^{-1}(x)\cap X$. 
We claim that there is an integer $2\leq l\leq 8s$ such
that if $b$ is a point in $\pi^{-1}(x)\cap X$ we have either
$\|a-b\|\leq (l-2)A\gamma^{1/s}t^{Nm/s}$ or
$\|a-b\|\geq (l+2)A\gamma^{1/s}t^{Nm/s}$. 
Indeed, $\pi^{-1}(x)\cap X$ contains only $s$ points and 
one can find in $\pi^{-1}(x)$ more than $s$ disjoint rings of the form
$$\Big\{w\in\pi^{-1}(x),\quad (l-2)A\gamma^{1/s}t^{Nm/s}\leq \|a-w\|< (l+2)A\gamma^{1/s}t^{Nm/s}\Big\}.$$
Let $B'$ be the
ball of center $x'$ and of radius
$lA\gamma^{1/s}t^{Nm/s}$ in $\C^k$ and $\partial B'$ its boundary. We
have $\dist(\pi^{-1}(x)\cap X, B\times\partial B')>A\gamma^{1/s}t^{Nm/s}$. The first
estimate on $\dist(\widetilde z,\pi_1^{-1}(x)\cap X)$ implies that
$X\cap B\times \partial B'=\varnothing$. So, $\pi$ is
proper on $X\cap B\times B'$ and it defines a ramified covering over
$B$. 

We show that the degree of this 
covering is at most equal to $m$. Otherwise, for some
multi-index $I$, 
$h_I(x)$ is  of order
$\gamma^{1/s}t^{Nm/s}\ll A^{-1}t^{N/s}=A^{-1}t^M$ because $\gamma$ is
small. This contradicts Lemma
\ref{lemma_technique} and the property that $\dist(x,Z)\geq t$. 
So, each connected component of $X\cap\pi^{-1}(B)$ is a covering of
degree $\leq m$ over $B$. Proposition \ref{prop_loj_1} implies the result.
\hfill $\square$ 

\medskip

The following corollary is deduced from Proposition \ref{prop_loj_2}. It is
enough to consider the restriction of the graph of $f$ to suitable
charts in $\P^k\times\P^k$ and then apply Proposition \ref{prop_loj_2} to each chart.

\begin{corollary} \label{cor_loj}
There is an integer $N\geq 1$ and a constant $A_4\geq 1$ such that if $0<t<1$ is a
constant and if $x,y$ are two points in
$\P^k$ with $\dist(x,E_\lambda)>t$ and
$\dist(y,E_\lambda)>t$, then we can write 
$$f^{-1}(x)=\{x_1,\ldots, x_{d^k}\}\quad \mbox{and}\quad 
f^{-1}(y)=\{y_1,\ldots, y_{d^k}\}$$
with $\dist(x_i,y_i)\leq A_4t^{-N}\dist(x,y)^{1/\delta_0}$.
\end{corollary}

We will use the following lemma in the proof of Theorem
\ref{th_equi_2}. 

\begin{lemma} \label{lemma_excep_bis}
There is a constant $A_5\geq 1$ and integers $L\geq 1$,  $n_1\geq 1$ such that if $a$ is a
  point out of $\Ec$, we can find a point $b\in f^{-n_1}(a)$ out of
  $E_\lambda$ such that 
$$\dist(b,E_\lambda)\geq A_5^{-1}\dist(a,\Ec)^L.$$
\end{lemma}
\proof
Let $F_n$ denote the set of points $a$ such that $f^{-n}(a)\subset
E_\lambda$.  
Since $E_\lambda$ is invariant, the sequence of analytic sets $F_n$ is
decreasing. By Lemma \ref{lemma_excep}, we
have $\cap F_n\subset \Ec$. So, $F_{n_1}\subset \Ec$ for $n_1$ large enough, that
is, $f^{-n_1}(a)\not\subset E_\lambda$ if $a\not\in\Ec$. 

Let $X$ denote the set of points $P=(a,b_1,\ldots,b_{d^{kn_1}})$ in
$\P^k\times (\P^k)^{d^{kn_1}}$ such that
$f^{-n_1}(a)=\{b_1,\ldots,b_{d^{kn_1}}\}$. Define $Y:=\P^k\times
E_\lambda^{d^{kn_1}}$ and $Z:=\Ec\times (\P^k)^{d^{kn_1}}$. We have
$X\cap Y\subset Z$. An inequality of Lojasiewicz
\cite[p.14 and p.62]{Malgrange} implies that 
$\dist(P,Y)\geq c\dist(P,X\cap Y)^L$ for some constants $c>0$ and
$L\geq 1$. We refer
$$d^{kn_1}\max_i\dist(b_i,E_\lambda)\geq \dist(P,Y)\geq
c\dist(P,Z)^L=c\dist(a,\Ec)^L.$$
The lemma follows.  
\endproof

\section{Equidistribution}

In this section, we prove Theorems \ref{th_equi_1} and \ref{th_equi_2}.
The theory of interpolation between Banach spaces allows to assume
that $\alpha=2$. Indeed, if $L:\Cc^0(\P^k)\rightarrow\R$ is a continuous
linear operator, then $L:\Cc^\alpha(\P^k)\rightarrow \R$ is also a
continuous linear operator and we have
$$\|L\|_{\Cc^\alpha}\leq
A\|L\|_{\Cc^0}^{1-\alpha/2}\|L\|_{\Cc^2}^{\alpha/2}$$
for some constant $A>0$ depending on $\alpha$ but independent of $L$,
see \cite{Triebel}.

We now prove Theorem \ref{th_equi_1}. Let $\Fc$ be the family of d.s.h. functions $\phi$ such that
$\langle\mu,\phi\rangle=0$ and $\|\phi\|_\mu\leq 1$.
Recall that the operator $f_*$ acts continuously on
d.s.h. functions. If $\phi$ is d.s.h., then 
$f_*(\phi)$ is defined by
$$f_*(\phi)(x):=\sum_{y\in f^{-1}(x)} \phi(y),$$
where the points in $f^{-1}(x)$ are counted with multiplicity. 
The following lemma is well-known and is easy to prove, see
\cite{DinhSibony6, DinhSibony4}.

\begin{lemma} \label{lemma_operator_lambda}
The operator $\Lambda:=d^{1-k}f_*$ acts
on d.s.h. functions and preserves the family $\Fc$.
\end{lemma}

We can assume that the function $\varphi$ in Theorems \ref{th_equi_1}
and \ref{th_equi_2} belongs to $\Fc$ and satisfies 
$\|\varphi\|_{\Cc^2}\leq 1$.
Let $A_i$, $\alpha$, $\eta$ and $N$ be the constants in Proposition \ref{prop_horm}, Lemmas
\ref{lemma_reg}, \ref{lemma_dist}, \ref{lemma_tub}, \ref{lemma_excep_bis} and 
 Corollary \ref{cor_loj}. 
Define 
$$l:=1+\log^+{1\over \dist(a,E_\lambda)}$$ 
where $a$ is a point out
of $E_\lambda$. 
Fix a  $\delta>1$ such that
$20k^2\delta_0<\delta<d\lambda^{-1}$.
Choose a constant $M>1$
large enough.
This choice depends on $\delta$. 
In what follows, the symbols
$\lesssim$ and $\gtrsim$ mean inequalities up to a multiplicative constant which
is independent of $\varphi$, $n$, $i$ and the variable $t$ that we
will use latter.

Fix an integer $n\geq 1$. For $0\leq i\leq n$, define
$\theta_i:=e^{-Ml\delta^in}$.
Define also $\varphi_0:=\varphi$ and write for $1\leq i\leq n$
$$\Lambda(\varphi_{i-1})=c_i+\varphi_i+\psi_i.$$
Here, $c_i+\varphi_i$ is the $\theta_i$-regularized function of
$\Lambda(\varphi_{i-1})$. The constant $c_i$ is chosen so that
$\langle \mu,\varphi_i\rangle=0$. Recall that we assumed that
$\langle\mu,\varphi\rangle=0$. We need the regularization because
$\Lambda(\varphi)$ is not $\Cc^1$ and we will apply Proposition \ref{prop_horm}.

\begin{lemma} \label{lemma_est_1}
The functions $\varphi_i$ and $\Lambda(\varphi_i)$ belong to $\Fc$. We
also have $c_i=-\langle\mu,\psi_i\rangle$ and
$$\|\nabla\varphi_i\|_\infty\lesssim \theta_i^{-7k^2},\quad
\|\varphi_i\|_\infty\leq (3d)^i \quad \mbox{and} \quad \|\psi_i\|_\infty\leq (3d)^i$$
for $1\leq i\leq n$.  
\end{lemma}
\proof
By Lemmas \ref{lemma_reg_dsh} and \ref{lemma_operator_lambda}, a
simple induction implies that $\varphi_i$ and $\Lambda(\varphi_i)$
belong to $\Fc$.
Hence, $\langle\mu,c_i+\psi_i\rangle=0$. It follows that $c_i=-\langle\mu,\psi_i\rangle$.
We also deduce that $\Lambda(\varphi_{i-1})$ has a bounded $L^1$-norm and by Lemma
\ref{lemma_reg}, $\|c_i+\varphi_i\|_{\Cc^1}\lesssim
\theta_i^{-7k^2}$. Therefore, 
$\|\nabla\varphi_i\|_\infty\lesssim \theta_i^{-7k^2}$. 

The last two inequalities are also proved by induction. Assume that
$|\varphi_{i-1}|\leq (3d)^{i-1}$ for some $1\leq i\leq n$. By definition of $\Lambda$, we have
$|\Lambda(\varphi_{i-1})|\leq 3^{i-1}d^i$. It follows from the definition of
$\theta$-regularization that $|c_i+\varphi_i|\leq 3^{i-1}d^i$ and then
$|\psi_i|\leq {2\over 3}(3d)^i\leq (3d)^i$. Since
$c_i=-\langle\mu,\psi_i\rangle$, we obtain $|c_i|\leq {2\over
  3}(3d)^i\leq (3d)^i$. Finally, we have 
$$|\varphi_i|\leq
|c_i+\varphi_i|+|c_i|\leq 3^{i-1}d^i+{2\over 3}
(3d)^i=(3d)^i.$$
This completes the proof.  
\endproof

We also have the following estimates where the constant $\eta$ has
been introduced in Lemma \ref{lemma_reg}.

\begin{lemma} \label{lemma_est_2}
For any $1\leq i\leq n$ and $t$ such that $2\eta\theta_i<t<1$, we have
$|\psi_i|\lesssim t^{-N}\theta_{i-1}^{-7k^2}\theta_i^{1/\delta_0}$ out
  of the $t$-neighbourhood of $E_\lambda$ and 
  $|c_i|\lesssim d^{-n}$. 
\end{lemma}
\proof
We will apply Corollary \ref{cor_loj} and use the inequality
$$|\phi(x_i)-\phi(y_i)|\lesssim\|\nabla
\phi\|_\infty\dist(x_i,y_i).$$ 
If $x,y$ are two points such that $\dist(x,E_\lambda)\geq {t\over
  2}$ and $\dist(y,E_\lambda)\geq {t\over
  2}$, then we deduce from the definition of $\Lambda$ that 
$$|\Lambda\varphi_{i-1}(x)-\Lambda\varphi_{i-1}(y)|\lesssim t^{-N}
\|\nabla\varphi_{i-1}\|_\infty\dist(x,y)^{1/\delta_0}\lesssim t^{-N} \theta_{i-1}^{-7k^2}
\dist(x,y)^{1/\delta_0}.$$ 
When $\dist(x,E_\lambda)\geq t$ and $\dist(x,y)\leq \eta \theta_i$, we
have    $\dist(y,E_\lambda)\geq {t\over 2}$ and then
$$|\Lambda\varphi_{i-1}(x)-\Lambda\varphi_{i-1}(y)|\lesssim t^{-N}
\theta_{i-1}^{-7k^2}\theta_i^{1/\delta_0}.$$ 
Recall that
$\psi_i=\Lambda(\varphi_{i-1})-\Lambda(\varphi_{i-1})_{\theta_i}$.  
We then deduce from Lemma \ref{lemma_reg} that 
$$|\psi_i(x)|\lesssim
t^{-N}\theta_{i-1}^{-7k^2}\theta_i^{1/\delta_0}$$ 
for $x$ out  of the $t$-neighbourhood of $E_\lambda$. 
This gives the first estimate in the lemma.

For the second one, choose $t:=\theta_i^{1/(2N\delta_0)}$. Let $V_t$ denote the
$t$-neighbourhood of $E_\lambda$. Using Lemmas \ref{lemma_tub} and \ref{lemma_est_1}, we obtain
\begin{eqnarray*}
|c_i| & \leq &  \int_{V_t} |\psi_i|d\mu + \int_{\P^k\setminus V_t}
|\psi_i| d\mu \lesssim (3d)^i\mu(V_t)+
t^{-N}\theta_{i-1}^{-7k^2}\theta_i^{1/\delta_0} \\
& \lesssim &  (3d)^i\theta_i^{\beta/(2N\delta_0)}+
\theta_{i-1}^{-7k^2}\theta_i^{1/(2\delta_0)}\lesssim d^{-n}.
\end{eqnarray*}
Here, we use the choice of $\theta_i$ and $\delta$.
\endproof

\noindent
{\bf End of the proof of Theorem \ref{th_equi_1}.}
Define $\nu_n:=d^{-kn} (f^n)^*(\delta_a)$. We need to show that
$|\langle \nu_n,\varphi\rangle|\leq A l\lambda^{-n}$
for some constant $A>0$. We have
\begin{eqnarray*}
\langle \nu_n,\varphi\rangle & = &  \langle d^{-k}
f^*(\nu_{n-1}),\varphi_0\rangle = \langle d^{-1}
\nu_{n-1},\Lambda(\varphi_0)\rangle\\
& = & d^{-1}c_1+ d^{-1}\langle
\nu_{n-1},\psi_1\rangle+ d^{-1}\langle
\nu_{n-1},\varphi_1\rangle.
\end{eqnarray*}
Using the identity $\Lambda(\varphi_{i-1})=c_i+\varphi_i+\psi_i$ and a
simple induction, we obtain
$$\langle \nu_n,\varphi\rangle = d^{-1}c_1+\cdots+d^{-n}c_n+ d^{-1}\langle\nu_{n-1},\psi_1\rangle+\cdots+
d^{-n} \langle \nu_0,\psi_n\rangle
+d^{-n}\langle\nu_0,\varphi_n\rangle.$$
By Lemma \ref{lemma_dist}, the distance between $\supp(\nu_i)$ and $E_\lambda$
is larger than $A_2^{-n}e^{-l}$. So, Lemma \ref{lemma_est_2} applied to
$t:={1\over 2}A_2^{-n}e^{-l}$ implies that 
$$|\langle \nu_{n-i},\psi_i\rangle|\lesssim A_2^{Nn} e^{Nl}
\theta_{i-1}^{-7k^2}\theta_i^{1/\delta_0}\lesssim d^{-n}.$$
We also have by Proposition \ref{prop_horm} and Lemma \ref{lemma_est_1} that
$$d^{-n}|\langle \nu_0,\varphi_n\rangle|\lesssim d^{-n}(1+\log^+\|\nabla\varphi_n\|)\lesssim
d^{-n}|\log\theta_n|\lesssim l\lambda^{-n}.$$
The previous estimates, together with the estimate on $c_i$ in Lemma
\ref{lemma_est_2}, imply the desired estimate $|\langle \nu_n,\varphi\rangle|\lesssim
l\lambda^{-n}$. 
\hfill $\square$

\

\noindent
{\bf Proof of Theorem \ref{th_equi_2}.} We can assume that $n$ is
even. We will apply Theorem
\ref{th_equi_1} and Lemma \ref{lemma_excep_bis}. In order to simplify
the notation, we can assume that $n_1=1$ in  Lemma
\ref{lemma_excep_bis} and use the above notations. 
So, for $a\not\in\Ec$, we have $|f^{-1}(a)\cap E_\lambda|\leq d^k-1$,
where we count the multiplicities of points.
Observe that if $x\not\in E_\lambda$ then $f^{-1}(x)\cap
E_\lambda=\varnothing$ since $E_\lambda$ is invariant. Therefore, 
$f^{-n/2}(a)\cap E_\lambda$ contains at
most $(d^k-1)^{n/2}$ points. 

Lemma \ref{lemma_excep_bis},
together with Lemma \ref{lemma_dist}, gives a more precise property. In fact,
there are at least  $d^{kn/2}-(d^k-1)^{n/2}$ points $b\in
f^{-n/2}(a)\setminus E_\lambda$, counted with
multiplicity, such that 
$$\dist(b,E_\lambda)\geq A_2^{-nL/2}
A_5^{-1}\dist(a,\Ec)^L\gtrsim A_2^{-nL/2} e^{-l'L},$$
where 
$$l':=1+\log^+{1\over \dist(a,\Ec)}\cdot$$  
Let $\Sigma$ denote the set of these points.
We have $d^{kn/2}-(d^k-1)^{n/2}\leq |\Sigma|\leq d^{kn/2}$ and 
$$d^{-kn}(f^n)^*(\delta_a)=d^{-kn}\sum_{b\in\Sigma}
(f^{n/2})^*(\delta_b) + d^{-kn}\sum_{b\in f^{-n/2}(a)\setminus\Sigma}
(f^{n/2})^*(\delta_b).$$
If $\varphi$ is as above, we deduce from Theorem
\ref{th_equi_1} applied to points $b\in\Sigma$ that
$$|\langle (f^{n/2})^*(\delta_b),\varphi\rangle|\lesssim d^{kn/2}
\Big[1+\log^+{1\over\dist(b,E_\lambda)}\Big]\lambda^{-n}\lesssim d^{kn/2}(n+l')\lambda^{-n}.$$
Hence, since $|f^{-n/2}(a)\setminus\Sigma|\leq (d^k-1)^{n/2}$ and
$|\varphi|\leq 1$, we obtain
\begin{eqnarray*}
|\langle  d^{-kn}(f^n)^*(\delta_a),\varphi\rangle| & \lesssim & 
d^{-kn}|\Sigma|d^{kn/2} (n+l') \lambda^{-n} +
d^{-kn} d^{kn/2}|f^{-n/2}(a)\setminus\Sigma|\\
& \lesssim & (n+l')\lambda^{-n} +(1-d^{-k})^{n/2}.
\end{eqnarray*}
It is enough to replace $\lambda$ with some smaller constant in order
to obtain the desired estimate $|\langle
d^{-kn}(f^n)^*(\delta_a),\varphi\rangle|\lesssim l'\lambda^{-n}$.
\hfill $\square$


T.-C. Dinh, UPMC Univ Paris 06, UMR 7586, Institut de
Math{\'e}matiques de Jussieu, F-75005 Paris, France. {\tt
  dinh@math.jussieu.fr}, {\tt http://www.math.jussieu.fr/$\sim$dinh} 

\

\noindent
N. Sibony,
Universit{\'e} Paris-Sud, Math{\'e}matique - B{\^a}timent 425, 91405
Orsay, France. {\tt nessim.sibony@math.u-psud.fr}

\end{document}